\def\version{0.64}
\def\journal{
ART
}
\def\titlep{An invariant of states on Cuntz algebras}
\documentclass[11pt]{article}
\usepackage{graphicx,ifthen,amssymb,amsmath,epic,eepic,color}
\usepackage{latexsym}

\newcommand{\qed}{\hbox{\rule[-2pt]{3pt}{6pt}}}
\newcommand{\qedh}{\hfill\qed \\}

\font\germ=eufm10 at12pt

\def\goth#1{\hbox{\germ#1}}

\newcommand{\vv}{\vspace{.3in}}

\newcommand{\vep}{\varepsilon}

\newtheorem{Thm}{Theorem}[section]

\newtheorem{rem}[Thm]{Remark}

\newtheorem{ex}[Thm]{Example}
\newtheorem{defi}[Thm]{Definition}
\newtheorem{lem}[Thm]{Lemma}

\newtheorem{prop}[Thm]{Proposition}

\newtheorem{cor}[Thm]{Corollary}
\newtheorem{fact}[Thm]{Fact}

\newcommand{\ww}{\vv\noindent}

\newcommand{\kn}{\Large\bf
$K\hspace{-.4cm} N$
\Large\bf\vv }

\def\cal#1{\mathcal #1}
\def\con{{\cal O}_{n}}

\def\edot{=1,\ldots,n}
\def\pr{{\it Proof.}\quad}

\def\evp{eventually periodic}

\def\co#1{{\cal O}_{#1}}
\def\ltn{\ell^{2}}
\def\ltno{\ell^{2}_1}

\addtocounter{footnote}{1}
\def\sftt#1{
\setcounter{equation}{0}
\addtocounter{footnote}{1}
\section{#1}
}

\def\ssft#1{\subsection{#1}}
\def\sssft#1{\subsubsection{#1}}

\begin{document}
%
%
\def\autherp{Katsunori Kawamura}
\def\emailp{e-mail: 
kawamurakk3@gmail.com
}
\def\addressp{{\small {\it College of Science and Engineering, 
Ritsumeikan University,}}\\
{\small {\it 1-1-1 Noji Higashi, Kusatsu, Shiga 525-8577, Japan}}
}

\def\scm#1{S({\Bbb C}^{N})^{\otimes #1}}
\def\mqb{\{(M_{i},q_{i},B_{i})\}_{i=1}^{N}}
\newcommand{\mline}{\noindent
\thicklines
\setlength{\unitlength}{.1mm}
\begin{picture}(1000,5)
\put(0,0){\line(1,0){1250}}
\end{picture}
\par
 }
\def\ptimes{\otimes_{\varphi}}
\def\delp{\Delta_{\varphi}}
\def\delf{\Delta_{{\Bbb F}}}
\def\delps{\Delta_{\varphi^{*}}}
\def\gamp{\Gamma_{\varphi}}
\def\gamps{\Gamma_{\varphi^{*}}}
\def\sem{\textsf{M}}
\def\hdelp{\hat{\Delta}_{\varphi}}
\def\tilco#1{\tilde{\co{#1}}}
\def\ba{\mbox{\boldmath$a$}}
\def\bb{\mbox{\boldmath$b$}}
\def\bc{\mbox{\boldmath$c$}}
\def\bd{\mbox{\boldmath$d$}}
\def\be{\mbox{\boldmath$e$}}
\def\bk{\mbox{\boldmath$k$}}
\def\bp{\mbox{\boldmath$p$}}
\def\bq{\mbox{\boldmath$q$}}
\def\bu{\mbox{\boldmath$u$}}
\def\bv{\mbox{\boldmath$v$}}
\def\bw{\mbox{\boldmath$w$}}
\def\bx{\mbox{\boldmath$x$}}
\def\by{\mbox{\boldmath$y$}}
\def\bz{\mbox{\boldmath$z$}}

\def\titlepage{

\noindent
{\bf 
\noindent
\thicklines
\setlength{\unitlength}{.1mm}
\begin{picture}(1000,0)(0,-100)
\put(0,0){\kn \knn\, for \journal\, Ver.\,\version}
\put(0,-50){\today,\quad {\rm file:}
 {\rm {\small \textsf{tit01.txt,\, J1.tex}}}}
\end{picture}
}
\vspace{-.5cm}
\quad\\
{\small 
\footnote{
\begin{minipage}[t]{6in}
directory: \textsf{\fileplace}, \\
file: \textsf{\incfile},\, from \startdate
\end{minipage}
}}
\quad\\
\framebox{
\begin{tabular}{ll}
\textsf{Title:} &
\begin{minipage}[t]{4in}
\titlep
\end{minipage}
\\
\textsf{Author:} &\autherp
\end{tabular}
}
{\footnotesize	
\tableofcontents }
}

\def\pdf#1{{\rm PDF}_{#1}}
\def\tilco#1{\tilde{\co{#1}}}
\def\ndm#1{{\bf M}_{#1}(\{0,1\})}
\def\cdm#1{{\cal M}_{#1}(\{0,1\})}
\def\tndm#1{\tilde{{\bf M}}_{#1}(\{0,1\})}
%
\def\openone{\mbox{{\rm 1\hspace{-1mm}l}}}
\def\goh{{\goth h}}
\def\gprs{geometric progression state}
\def\Gpr{Geometric progression\ }
\def\qmatrix#1{
\left[
\begin{matrix}
#1
\end{matrix}
\right]}
\def\csis{cyclic $s_i^*$-invariant subspace}
\def\tr{{\rm Tr}}
\def\ip#1#2{\langle #1|#2\rangle}
\def\csa{{\Bbb C}\langle s_1^*,\ldots, s_n^*\rangle}
\def\cfnp{{\Bbb C}{\Bbb F}_n^+}
\newcommand{\bs}{\boldsymbol{\sigma}}
\newcommand{\Hom}{\mathop{\mathrm{Hom}}\nolimits}
\def\lspan#1{\mathop{\mathrm{Lin}}\nolimits\langle\{#1\}\rangle}
\def\Ind{\mathop{\mathrm{Ind}}\nolimits}
\def\cdim{{\rm c}\!\dim}
\def\End{\mathop{\mathrm{End}}\nolimits}
\def\Rep{\mathop{\mathrm{Rep}}\nolimits}
\def\Irr{\mathop{\mathrm{Irr}}\nolimits}
\def\Erg{\mathop{\mathrm{Erg}}\nolimits}
\def\bh{{\cal B}({\cal H})}
\def\endbh{\End{\cal B}({\cal H})}
\def\endnbh{\End_n{\cal B}({\cal H})}
%
%
%
\setcounter{section}{0}
\setcounter{footnote}{0}
\setcounter{page}{1}
\pagestyle{plain}

%
%
\title{\titlep}
\author{\autherp\thanks{\emailp}
\\
\addressp}
\date{}
\maketitle
%
%
\begin{abstract}
For an arbitrary state $\omega$ on a Cuntz algebra,
we define a number $1\leq \kappa(\omega)\leq \infty$
such that 
if the GNS representations of $\omega$ and $\omega'$
are unitarily equivalent, then $\kappa(\omega)=\kappa(\omega')$.
By using $\kappa$,
we define minimal states and
it is shown that
the classification problem of states is reduced to that of minimal states.
By using results of Dutkay, Haussermann, and Jorgensen,
we give a sufficient condition of the minimality of a state.
Properties of $\kappa$ and examples are shown.
As an application, a new invariant of a certain class of endomorphisms of 
${\cal B}({\cal H})$ is given.
\end{abstract}


%
%
\sftt{Introduction}
\label{section:first}
The most different aspect in operator algebra from other mathematics
is the treatment of non-type I C$^*$-algebras \cite{Glimm}.
By definition, 
a non-type I C$^*$-algebra is characterized by its representations.
Hence the study of representations of non-type I C$^*$-algebras
is a core component of operator algebra.
For example,  Cuntz algebras are non-type I.
The aim of this paper is to classify states on Cuntz algebras
by using a new invariant.
In this section, we introduce the invariant and show its properties.
In $\S$ \ref{subsection:firsttwo}, we will state our main results.
In $\S$ \ref{subsection:firstthree}, 
the significance and advantages of the new invariant 
will be explained.
%
%
\ssft{Invariant}
\label{subsection:firstone}
%
%
\sssft{Definition}
\label{subsubsection:firstoneone}
For $2\leq n\leq \infty$,
let $\con$ denote the {\it Cuntz algebra} 
with Cuntz generators $s_1,\ldots,s_n$ \cite{C}, 
that is, 
$\con$ is a C$^{*}$-algebra 
which is universally generated by a (finite or infinite) sequence
$s_{1},\ldots,s_n$ satisfying
$s_{i}^{*}s_{j}=\delta_{ij}I$ for $i,j\edot$ and
%
%
\begin{equation}
\label{eqn:cuntzdef}
\sum_{i=1}^{n}s_{i}s_{i}^{*}=I\mbox{ when } n<\infty,\quad
\sum_{i=1}^{k}s_{i}s_{i}^{*}\leq I,\quad k= 1,2,\ldots 
\mbox{ when }n = \infty
\end{equation}
where $I$ denotes the unit of $\con$.
The Cuntz algebra $\con$ is an infinite dimensional, 
noncommutative simple C$^{*}$-algebra with unit.

Let ${\cal S}(\con)$ denote
the set of all states on $\con$.
For $\omega,\omega'\in {\cal S}(\con)$,
we write $\omega\sim\omega'$ when
their Gel'fand-Naimark-Segal (=GNS)  
 representations are unitarily equivalent.
The problem is to classify elements in ${\cal S}(\con)$ by 
the equivalence relation $\sim$.
For $\omega\in {\cal S}(\con)$
with GNS representation $({\cal H},\pi,\Omega)$,
define the nonzero closed subspace ${\cal K}(\omega)$ of ${\cal H}$ 
(\cite{BJKW,BJP,DHJ}) by
%
%
\begin{equation}
\label{eqn:komega}
{\cal K}(\omega):=\overline{\lspan{\pi(s_J)^*\Omega:J\in {\cal I}_n}}
\end{equation}
where ${\cal I}_n:=\bigcup_{l\geq 0}\{1,\ldots,n\}^l$,
$\{1,\ldots,n\}^0:=\{\emptyset\}$,
$s_J:=s_{j_1}\cdots s_{j_l}$ for $J=(j_1,\ldots,j_l)$,
and $s_{\emptyset}:=I$.
When $n=\infty$,
replace $\{1,\ldots,n\}^l$ with $\{1,2,\ldots\}^l$.
Define $\cdim \omega$ and $\kappa(\omega)$ by 
%
%
\begin{equation}
\label{eqn:cdimkappa}
\cdim\omega:=\dim {\cal K}(\omega),\quad
\kappa(\omega):=\min\{\cdim\omega':
\omega'\in {\cal S}(\con),\, \omega'\sim\omega\}.
\end{equation}
A state $\omega$ on $\con$ is said to be {\it minimal}
if $\cdim\omega=\kappa(\omega)$.
By definition, the following hold immediately.
%
%
\begin{Thm}
\label{Thm:mainone}
\begin{enumerate}
\item
\label{Thm:mainoneone}
For any $\omega\in {\cal S}(\con)$,
there exists a minimal state $\omega'$ on $\con$
which is equivalent to $\omega$.
We call such $\omega'$ a {\it minimal model}  of $\omega$.
\item
\label{Thm:mainonetwo}
For $\omega,\omega'\in {\cal S}(\con)$,
if $\omega\sim\omega'$,
then $\kappa(\omega)=\kappa(\omega')$.
\end{enumerate}
\end{Thm}
%
%
\pr
(i)
Since $\{\cdim \omega':\omega'\sim\omega\}$
is a subset of $\{1,2,\ldots,\infty\}$,
it always has the smallest element
with respect to the standard linear ordering
where $\infty$ means the countably infinite cardinality.
Hence there always exists  a minimal state $\omega'$
which is equivalent to $\omega$.

\noindent
(ii)
Assume $\omega\sim\omega'$.
Then their minimal models are also equivalent.
By definitions of $\kappa(\omega)$
and $\kappa(\omega')$,
the statement holds.
\qedh

\noindent
From Theorem \ref{Thm:mainone},
the classification problem of (pure) states on $\con$
is reduced to that of minimal (pure) states on $\con$
with $\cdim =d$ for each number $1\leq d\leq\infty$.
Remark that $\kappa(\omega)$ is an invariant of 
$\omega$, but 
$\cdim \omega$ is not (Proposition \ref{prop:difference}).
For a given $1\leq d\leq \infty$,
there exist continuously many minimal pure states 
with $\kappa(\omega)=d$ (Theorem \ref{Thm:conti}).
A minimal model of a state is not unique in general
(Proposition \ref{prop:uniqueness}).
The symbol ``$\cdim\omega$" originates in 
our old terminology,
``the correlation dimension of $\omega$" 
(see Definition \ref{defi:proper}(i)).

For $2\leq n\leq \infty$,
we write $U(n)$ as the rank $n$ unitary group when $n<\infty$,
as the group of all unitaries on 
$\ltn:=\{(z_j): \sum_{j\geq 1}|z_j|^2=1\}$ when $n=\infty$.
We show properties of $\kappa$ 
with respect to the unitary group action as follows.
%
%
\begin{prop}
\label{prop:unitary}
{\rm (}$U(n)$ invariance{\rm )}
Let $\alpha$ denote the standard $U(n)$-action
on $\con$, that is,
$\alpha_g(s_i):=\sum_{j=1}^ng_{ji}s_j$
for $i=1,\ldots,n$ and $g=(g_{ij})\in U(n)$.
Let $\omega\in {\cal S}(\con)$.
\begin{enumerate}
\item
For any $g\in U(n)$,
$\cdim (\omega\circ \alpha_g)=\cdim\omega$.
\item
\label{prop:unitarytwo}
For any $g\in U(n)$,
$\kappa(\omega\circ  \alpha_g)=\kappa(\omega)$.
\end{enumerate}
\end{prop}
%
%
\pr
(i)
Let $({\cal H},\pi,\Omega)$ denote 
the GNS representation of $\omega$.
Since
$\omega=\ip{\Omega}{\pi(\cdot)\Omega}$,
we obtain $\omega\circ \alpha_g=\ip{\Omega}{\pi(\alpha_g(\cdot))
\Omega}$.
Hence we identify
the GNS representation of $\omega\circ \alpha_g$
with $({\cal H},\pi\circ \alpha_g,\Omega)$
(see 4.5.3 Proposition of \cite{KR1}).
Then ${\cal K}(\omega\circ \alpha_g)$
is spanned by
the set $\{\pi(\alpha_g(s_J))^*\Omega:J\in {\cal I}_n\}$.
This is contained in ${\cal K}(\omega)$
by the definition of $\alpha_g$.
From this,
${\cal K}(\omega\circ \alpha_g)\subset {\cal K}(\omega)$.
By replacing $(\omega,g)$ with
$(\omega\circ \alpha_g, g^*)$,
we obtain
${\cal K}(\omega)
={\cal K}((\omega\circ \alpha_g)\circ \alpha_{g^*})\subset 
{\cal K}(\omega\circ \alpha_g)$.
Hence the statement holds.

\noindent
(ii)
Remark that $\omega\sim \omega'$ if and only if 
$\omega\circ \alpha_g\sim \omega'\circ \alpha_g$.
From this and (i), the statement holds.
\qedh

\noindent
From Proposition \ref{prop:unitary}
and Theorem \ref{Thm:mainone},
$\kappa(\omega)$ can be regarded 
as an invariant of a $U(n)$-orbit 
in the set of all unitary equivalence classes of (pure) states on $\con$.
%
%
\begin{cor}
\label{cor:unitary}
For $\omega,\omega'\in {\cal S}(\con)$,
if $\omega'\sim \omega\circ \alpha_g$ for some $g\in U(n)$,
then
$\kappa(\omega)=\kappa(\omega')$.
\end{cor}
%
%
\pr
From Proposition \ref{prop:unitary}(ii) and 
Theorem \ref{Thm:mainone}(ii),
the statement holds.
\qedh

%
%
\sssft{Cuntz states as the case of $\kappa=1$}
\label{subsubsection:firstonetwo}
We review well-known results about Cuntz states
by using $\cdim$ and $\kappa$.
Let $({\Bbb C}^n)_1:=\{z\in {\Bbb C}^n:\|z\|=1\}$.
For any $z=(z_1,\ldots,z_n)\in ({\Bbb C}^n)_1$,
a state $\omega_z$ on $\con$ which satisfies
%
%
\begin{equation}
\label{eqn:first}
\omega_z(s_j)=\overline{z_j}\quad\mbox{for all }j=1,\ldots,n,
\end{equation}
exists uniquely and is pure,
where $\overline{z_j}$ denotes the complex conjugate of $z_j$.
When $n=\infty$, replace ${\Bbb C}^n$ with $\ltn$.
The state $\omega_z$ is called the {\it Cuntz state} 
by $z$ \cite{BEGJ,BJ1997,BJ,BJP,BK,CE,Cuntz1980,GP08,GP15}.
For any $g\in U(n)$ and $z\in ({\Bbb C}^n)_1$,
$\omega_z\circ \alpha_{g^*}=\omega_{gz}$,
that is, the group $U(n)$ transitively acts on the set of 
all Cuntz states on $\con$.
%
%
\begin{Thm}
\label{Thm:Cuntzstateeq}
{\rm (\cite{GP08}, Appendix B)}
For $z,y\in ({\Bbb C}^n)_1$,
$\omega_z\sim\omega_y$ if and only if $z=y$.
\end{Thm}
From Theorem \ref{Thm:Cuntzstateeq},
$({\Bbb C}^n)_1$ is the complete set of invariants
of Cuntz states on $\con$.
%
%
\begin{fact}
\label{fact:seven}
Assume $2\leq n\leq \infty$.
\begin{enumerate}
\item
For $\omega\in {\cal S}(\con)$,
$\cdim\omega=1$ if and only if $\omega$ is a Cuntz state.
Especially, any Cuntz state is minimal.
\item
For $\omega\in {\cal S}(\con)$,
$\kappa(\omega)=1$ if and only if 
$\omega$ is equivalent to a Cuntz state.
\item
For $\omega,\omega'\in {\cal S}(\con)$,
if $\kappa(\omega)=1=\kappa(\omega')$,
then $\omega'\sim \omega\circ \alpha_g$ for some $g\in U(n)$.
\end{enumerate}
\end{fact}
%
%
\pr
(i)
Let $({\cal H},\pi,\Omega)$ denote the GNS representation 
of $\omega$.
We see that 
(\ref{eqn:first}) is equivalent 
that $\pi(s_j)^*\Omega=z_j\Omega$ for all $j$.
From this and the definition of $\cdim$, the statement holds.

\noindent
(ii) 
By definition,
$\kappa(\omega)=1$
if and only if $\omega$ is equivalent to $\omega'$ such that
$\cdim \omega'=1$.
This is equivalent to the statement that 
$\omega$ is equivalent to a Cuntz state from (i).

\noindent
(iii) 
From (ii),
there exist Cuntz states $\omega_1$ and $\omega_1'$ 
such that $\omega\sim \omega_1$ and $\omega'\sim \omega_1'$.
Since
$\omega_1\circ \alpha_g=\omega_1'$ for some $g\in U(n)$,
the statement holds.
\qedh

\noindent
By combining Fact \ref{fact:seven}(ii) 
and Theorem \ref{Thm:Cuntzstateeq},
the case of $\kappa=1$ is completely classified.
Remark that $\kappa=1$ implies the purity of a state automatically
because any Cuntz state is pure.
Fact \ref{fact:seven}(iii) does not hold
when $\kappa(\omega)=\kappa(\omega')\geq 2$
(Example \ref{ex:first}).

Remark that any Cuntz state is completely defined
by only a parameter $z\in ({\Bbb C}^n)_1$.
This is stated as the ``uniqueness" of $\omega_z$
in (\ref{eqn:first}).
In other words, it is not necessary to define
the value $\omega_z(s_Js_K^*)$ for all $J,K\in {\cal I}_n$.
Thanks to the uniqueness,
one can describe Cuntz states very concisely.
This type uniqueness holds for various 
other states in $\S$ \ref{section:third}.

%
%
\ssft{Main theorems}
\label{subsection:firsttwo}
We state our main theorems in this subsection.
Since their proofs require some lemmas,
we will prove theorems in $\S$ \ref{subsection:secondtwo}.
%
%
\sssft{Minimality of a state}
\label{subsubsection:firsttwoone}
Let $\cdim\omega$ and 
$\kappa(\omega)$ be as in (\ref{eqn:cdimkappa}).
In order to make use of our new invariant $\kappa$, 
we must be able to compute $\kappa(\omega)$.
If we know that $\omega$ is minimal,
then $\cdim \omega=\kappa(\omega)$.
Since the computation of $\cdim \omega$ is easier than
that of $\kappa(\omega)$,
the determination of its minimality makes sense.

Let ${\cal I}_n$ be as in (\ref{eqn:komega}).
Define
%
%
\begin{equation}
\label{eqn:conplus}
\con^+:=\overline{\lspan{s_J:J\in {\cal I}_n,\, J\ne \emptyset}}
\subset \con.
\end{equation}
Then $\con^+$ is a nonunital,
non-selfadjoint subalgebra of $\con$.
We obtain a sufficient condition
that a given state is minimal.
%
%
\begin{Thm}
\label{Thm:minimal}
For $\omega\in {\cal S}(\con)$,
if there exists an isometry $u$ in $\con^+$ 
{\rm (}that is, $u^*u=I${\rm )}
such that $\omega(u)=1$,
then $\omega$ is minimal.
\end{Thm}

\noindent
By using Theorem \ref{Thm:minimal},
we will show examples of minimal state in 
$\S$ \ref{subsection:thirdone}.
For a state,  its minimality is neither necessary
nor sufficient for its purity
in general (Proposition \ref{prop:differenceb}).

%
%
\sssft{Properly infinite correlation of a state}
\label{subsubsection:firsttwotwo}
%
%
\begin{defi}
\label{defi:proper}
\begin{enumerate}
\item {\rm (\cite{BJ1997})}
A state $\omega$ on $\con$
is said to be {\it infinitely correlated}
if $\cdim \omega=\infty$.
Otherwise,
$\omega$ 
is said to be {\it finitely correlated}.
\item
%
A state $\omega$ on $\con$
is said to be {\it properly infinitely correlated}
if $\kappa(\omega)=\infty$.
Otherwise,
$\omega$ 
is said to be {\it essentially finitely correlated}.
\end{enumerate}
\end{defi}
By definition,
a state is either properly infinitely correlated
or essentially finitely correlated.
Any finitely correlated state is essentially finitely
correlated,
but the converse is not true
(Proposition \ref{prop:difference}).
If $\omega$ is pure and finitely correlated,
then any vector state of 
the GNS representation space of $\omega$
is essentially finitely correlated.
From Theorem \ref{Thm:mainone}(ii), the following holds.
%
%
\begin{fact}
\label{fact:ess}
For $\omega,\omega'\in {\cal S}(\con)$,
assume $\omega\sim\omega'$.
If $\omega$ is properly infinitely correlated,
then so is $\omega'$.
Otherwise,  $\omega'$ is essentially finitely correlated.
\end{fact}

We give a sufficient condition such that
a given state is properly infinitely correlated.
%
%
\begin{Thm}
\label{Thm:inclusion}
Let $\con^+$ be as in  {\rm (\ref{eqn:conplus})}.
For $\omega\in {\cal S}(\con)$,
assume that there exists a sequence $(a_i)_{i\geq 1}$
of isometries in $\con^+$ which satisfies
%
%
\begin{equation}
\label{eqn:aoneatwo}
\omega(a_1\cdots a_l\, a_k^*\cdots a_1^*)
=\delta_{lk}\quad \mbox{ for all }l,k\geq 1.
\end{equation}
Then $\omega$ is properly infinitely correlated.
\end{Thm}
By using Theorem \ref{Thm:inclusion},
we will show examples of properly infinitely correlated
state in $\S$ \ref{subsection:thirdtwo}.

%
%
\ssft{Summary of results}
\label{subsection:firstthree}
We summarize the significance and advantages of $\kappa$.
\begin{enumerate}
\item {\bf Refinement of definitions}:
According to \cite{BJ1997},
there exist two classes of states on $\con$,
that is, finitely correlated states and infinitely correlated states
(Definition \ref{defi:proper}(i)).
From Proposition \ref{prop:difference},
it becomes clear that this classification is incompatible 
with the unitary equivalence of states.
For example, Figure 1 in \cite{GP08} is exceedingly inappropriate.
Instead of these notions,
essentially finitely/properly infinitely correlated states
are established by using $\kappa$ (Fact \ref{fact:ess}).
\item {\bf New classification method}:
As a classification theory of representations of C$^*$-algebras,
the Murray-von Neumann-Connes classification 
\cite{MN01, Connes} is well known,
that is, for a given factor representation,
its type is determined by 
the type of the von Neumann algebra generated by its image. 
Unfortunately, this classification is no use for
the classification of irreducible representations of $\con$
because the type of any irreducible representation of $\con$ 
is type I$_{\infty}$.
As a finer classification of irreducible representations of $\con$,
$\kappa$ is essentially new ($\S$ \ref{subsection:fourthone}).
\item {\bf Invariant for arbitrary states}:
Until now, 
there exist several small subclasses of states 
or representations of $\con$
which are often completely classified 
\cite{ABC,BC,BC2,BJ,Cuntz1980,DaPi2,DaPi3,DJ,TS11,GP08,GP15}.
They are often parameterized by their complete sets of invariants.
On the other hand,
$\kappa$ can be defined on the whole of states
(see also (\ref{eqn:widecon})).
\item
{\bf Reduction}:
The new notion ``minimal state" reduces the classification problem of 
states to that of minimal states from Theorem \ref{Thm:mainone}(i).
\item
{\bf Many examples}:
The existence of many examples firmly establishes 
that the theory of $\kappa$ is not vacuous.
In $\S$ \ref{subsubsection:firstonetwo} and $\S$ \ref{section:third},
examples are shown and their values of $\kappa$ are computed.
Especially,  we will show that the cardinality of the set of 
mutually inequivalent pure states with same invariant number 
are uncountable in Theorem \ref{Thm:conti}.
\item
{\bf Generalization of symbolic dynamical system}:
In Remark \ref{rem:dynamical},
we will show that $\kappa$ is a 
generalization of period length in the full one-sided shift.
This implies a naturality of $\kappa$.
\item
{\bf Applications}:
In $\S$ \ref{section:fourth},
we will show that $\kappa$
can be defined on both 
arbitrary irreducible representations of $\con$
and 
arbitrary ergodic endomorphisms of ${\cal B}({\cal H})$
as their invariants.
\end{enumerate}

The paper is organized as follows.
In $\S$ \ref{section:second},
we will prove Theorem \ref{Thm:minimal}
and Theorem \ref{Thm:inclusion}.
In $\S$ \ref{section:third},
we will show examples.
In $\S$ \ref{section:fourth},
we will show applications.

%
%
\sftt{Proofs of main theorems}
\label{section:second}
In this section,
we prove Theorem \ref{Thm:minimal}
and Theorem \ref{Thm:inclusion}.
For this purpose, we prove lemmas needed later.
%
%
\ssft{Dutkay-Haussermann-Jorgensen theory and its generalization}
\label{subsection:secondone}
We review a part of the work by
Dutkay, Haussermann, and Jorgensen (\cite{DHJ}, $\S$ 3.1) 
which shows a kind of structure theorem 
of a representation space of $\con$ in a general setting.
Our analysis is dependent on their results to a great extent.
We write ``a representation of $\con$"
to denote a unital $*$-representation of $\con$ in this paper.
%
%
\sssft{Dutkay-Haussermann-Jorgensen decomposition}
\label{subsubsection:secondoneone}
Let ${\cal K}(\omega)$ and ${\cal I}_n$ be as in 
(\ref{eqn:komega}).
%
%
\begin{defi}
\label{defi:csis}
Let $({\cal H},\pi)$ be a representation on $\con$ and $M$ a subspace.
\begin{enumerate}
\item 
{\rm (\cite{BJKW}, $\S$ 1)}
$M$ is said to be {\it $s_i^*$-invariant} 
if $\pi(s_i)^*M\subset M$ for all $i$.
\item
{\rm (\cite{DHJ}, Definition 2.4)}
$M$ is said to be {\it cyclic} if 
$\{\pi(s_Js_K^*)x: J,K\in {\cal I}_n,\,x\in M\}$
spans ${\cal H}$.
\end{enumerate}
\end{defi}

\noindent
Both $\{0\}$ and ${\cal H}$ are trivial $s_i^*$-invariant 
subspaces of ${\cal H}$.
Therefore any nonzero representation of $\con$ contains
a nonzero $s_i^*$-invariant subspace at least.
If $M$ contains a cyclic vector, then $M$ is cyclic.
If $({\cal H},\pi)$ is irreducible,
then any nonzero subspace of ${\cal H}$ is cyclic.
For any state $\omega$,
${\cal K}(\omega)$ in (\ref{eqn:komega})
is a closed \csis\ of the GNS representation space of $\omega$.
%
%
\begin{rem}
\label{rem:one}
{\rm
Let ${\cal R}_n\subset \con$ denote the algebra generated by 
$s_1^*,\ldots,s_n^*$ over ${\Bbb C}$, that is, 
%
%
\begin{equation}
\label{eqn:rn}
{\cal R}_n:={\Bbb C}\langle s_1^*,\ldots,s_n^* \rangle
\end{equation}
which consists of all noncommutative polynomials in 
$s_1^*,\ldots,s_n^*$ over ${\Bbb C}$.
We see that 
${\cal R}_n$ is the free algebra generated by 
$s_1^*,\ldots,s_n^*$ over ${\Bbb C}$
(\cite{Cohn}, $\S$ 6.2).
Remark that ${\cal R}_n$ is not a self-adjoint algebra
because $s_1,\ldots,s_n\not\in {\cal R}_n$.
Clearly, 
the standard terminology of 
``$s_i^*$-invariant subspace" is  just ``left ${\cal R}_n$-module."
We use the conventional word ``$s_i^*$-invariant" in this paper.
}
\end{rem}
%
%
\begin{Thm}{\rm (Dutkay-Haussermann-Jorgensen \cite{DHJ})}
\label{Thm:dhjone}
Let $({\cal H},\pi)$ be a representation of $\con$.
If $M$ is a closed \csis\ of ${\cal H}$,
then there exists a unique orthogonal decomposition of ${\cal H}$,
%
%
\begin{equation}
\label{eqn:dhjdeco}
{\cal H}=\bigoplus_{l\geq 0}{\cal H}_l
\end{equation}
such that
$\{\pi(s_J)v: 
J\in\{1,\ldots,n\}^k,\, v\in M\}$ 
spans $\bigoplus_{l=0}^k{\cal H}_l$
for all $k\geq 0$.
We call {\rm (\ref{eqn:dhjdeco})} the 
{\it Dutkay-Haussermann-Jorgensen (=DHJ) decomposition}
of $({\cal H},\pi)$ by $M$.
\end{Thm}
%
%
\pr
The existence is proved in $\S$ 3.1 of \cite{DHJ}.
The uniqueness holds from the properties of 
subspaces $\bigoplus_{l=0}^k{\cal H}_l$ for $k\geq 0$.
\qedh

\noindent
Let $\con^+$ be as in (\ref{eqn:conplus}).
By definition, ${\cal H}_0=M$ and 
$\pi(s_i)^*{\cal H}_l\subset {\cal H}_{l-1}$ 
and $\pi(s_i){\cal H}_l\subset {\cal H}_{l+1}$ 
for all $i$ when $l\geq 1$.
From this, any $x\in \con^+$ satisfies 
%
%
\begin{equation}
\label{eqn:pixstar}
\pi(x^*){\cal H}_l\, \subset\, \bigoplus_{k=0}^{l-1}{\cal H}_{k}
\quad(l\geq 1).
\end{equation}
Theorem \ref{Thm:dhjone} 
indicates that an essential part of a representation of $\con$
is its $s_i^*$-invariant subspace.
In Theorem \ref{Thm:dhjone},
$\pi(s_i){\cal H}_0$ is not a subspace of ${\cal H}_{1}$ 
in general (Proposition \ref{prop:dhjb}).
If $M={\cal H}$,
then ${\cal H}_l=\{0\}$ for all $l\geq 1$.

%
%
\sssft{Lemmas}
\label{subsubsection:secondonetwo}
The following are slight generalizations
of Theorem 3.2 and Lemma 3.4 in \cite{DHJ}.
%
%
\begin{lem}
\label{lem:appone}
Let $\con^+$ be as in 
{\rm (\ref{eqn:conplus})}.
Assume that $\con$ acts on a Hilbert space ${\cal H}$,
$M$ is a closed \csis\ of ${\cal H}$, and 
$a=(a_i)_{i\geq 1}$ is a sequence of isometries in $\con^+$.
Let $a[l]:=a_1\cdots a_l$ for $l\geq 1$.
\begin{enumerate}
\item
For any $\vep>0$ and $v\in {\cal H}$,
there exists $l_0\geq 1$
such that 
%
%
\begin{equation}
\label{eqn:pmai}
\|P_M a[l]^*v-a[l]^*v\|<\vep
\quad \mbox{for all $l\geq l_0$}
\end{equation}
where $P_M$ denotes the projection from ${\cal H}$
onto $M$.
\item
Define the projection 
$T:=\bigwedge_{l\geq 1} a[l]a[l]^*$ on ${\cal H}$.
If $\dim M<\infty$, then
$T{\cal H}\subset M$.
\end{enumerate}
\end{lem}
%
%
\pr
(i)
Let ${\cal H}=\bigoplus_{l\geq 0}{\cal H}_l$ 
denote the DHJ decomposition by $M$
 (Theorem \ref{Thm:dhjone}).
From (\ref{eqn:pixstar}),
$a_i^*{\cal H}_l\subset 
\bigoplus_{k=0}^{l-1}{\cal H}_k$ for all $i\geq 1$
when $l\geq 1$.
We can find $l_0\geq 1$ and vectors $v_1,v_2\in {\cal H}$ such that 
$v=v_1+v_2$,
$v_1\in \bigoplus_{l=0}^{l_0}{\cal H}_l$,
$v_2\in \bigoplus_{l>l_0}{\cal H}_l$,
and $\|v_2\|<\vep$.
Then
$\|v_2\|<\vep$ and $a[l]^*v_1\in {\cal H}_0=M$ 
for all $l\geq l_0$.
Hence $(P_M-I)a[l]^*v_1=0$ for all $l\geq l_0$.
From this,
%
%
\begin{equation}
\label{eqn:pmalv}
\begin{array}{rl}
P_Ma[l]^*v-a[l]^*v
=&
(P_M-I)a[l]^*v\\
=&
(P_M-I)a[l]^*(v_1+v_2)\\
=&
(P_M-I)a[l]^*v_2.
\end{array}
\end{equation}
Therefore
%
%
\begin{equation}
\label{eqn:pmaep}
\|P_M a[l]^*v-a[l]^*v\|
\leq \|P_M-I\| \  \|a[l]^*\| \ \|v_2\|
\leq \|v_2\|<\vep.
\end{equation}
\noindent
(ii)
Remark that 
$a[l]^*a[l]=I$ for all $l$
because $a[l]$ is a product of isometries.
Since $\{a[l]a[l]^*:l\geq 1\}$ is a decreasing sequence of 
projections on ${\cal H}$, $T$ is well defined.
It is sufficient to show the case of $T\ne 0$.
Assume $T\ne 0$.
We prove 
$(T{\cal H}\cap M)^{\perp}\cap T{\cal H}=\{0\}$.
By definitions of $T$ and $a[l]$,
we see that $a[l]^*$ is a unitary on $T{\cal H}$ for all $l\geq 1$.
Since 
$a[l]^*(T{\cal H}\cap M)\subset T{\cal H}\cap M$
and $\dim M<\infty$,
we obtain $a[l]^*(T{\cal H}\cap M)= T{\cal H}\cap M$.
This implies 
$a[l]^*\{(T{\cal H}\cap M)^{\perp}\cap T{\cal H}\}
= (T{\cal H}\cap M)^{\perp}\cap T{\cal H}$.
Hence $a[l]^*$ is also a unitary from 
$(T{\cal H}\cap M)^{\perp}\cap T{\cal H}$
onto $ (T{\cal H}\cap M)^{\perp}\cap T{\cal H}$.

Assume $v\in (T{\cal H}\cap M)^{\perp}\cap T{\cal H}$.
Then 
%
%
\begin{equation}
\label{eqn:alv}
a[l]^*v\in (T{\cal H}\cap M)^{\perp}\cap T{\cal H}\quad
(l\geq 1).
\end{equation}
From (i),
%
%
\begin{equation}
\label{eqn:convtwo}
\|a[l]^*v-P_Ma[l]^*v\|\to 0\quad\mbox{ when }l\to\infty.
\end{equation}
Since $\dim M<\infty$,
the sequence $\{P_Ma[l]^*v:l\geq 1\}$
has a convergent subsequence 
$\{P_Ma[l_i]^*v:i\geq 1\}$
in the compact subset 
$M':=\{x\in M:\|x\|\leq \|v\|\}$ of $M$.
Let
%
%
\begin{equation}
\label{eqn:vlim}
v_{\infty}:=\lim_{i\to\infty}P_Ma[l_i]^*v\in P_M{\cal H}=M.
\end{equation}
From this, (\ref{eqn:convtwo}), and (\ref{eqn:alv}),
we obtain
$v_{\infty}=\lim_{i\to\infty}a[l_i]^*v\in 
(M\cap T{\cal H})^{\perp}\cap T{\cal H}$.
From this and (\ref{eqn:vlim}),
%
%
\begin{equation}
\label{eqn:vinf}
v_{\infty}\in 
M\cap (M\cap T{\cal H})^{\perp}\cap T{\cal H}
=
(M\cap T{\cal H})^{\perp}\cap (M\cap T{\cal H})
=\{0\}.
\end{equation}
Hence $v_{\infty}=0$.
On the other hand,
since $a[l]^*$ is a unitary on 
$(M\cap T{\cal H})^{\perp}\cap T{\cal H}$,
we obtain
$0=\|v_{\infty}\|
=\lim _{i\to\infty}\|a[l_i]^*v\|=\|v\|$.
Hence $v=0$.
Therefore 
$(T{\cal H}\cap M)^{\perp}\cap T{\cal H}=\{0\}$.
\qedh

\noindent
In (\ref{eqn:pmaep}), if $M={\cal H}$, then $\|P_M-I\| =0$. 
Hence 
``$\|P_M-I\| \  \|a[l]^*\| \ \|v_2\| \leq \|v_2\|$"
can not be replaced with 
``$\|P_M-I\| \  \|a[l]^*\| \ \|v_2\| = \|v_2\|$" in general.
%
%
\begin{lem}
\label{lem:previous}
Assume that $\con$ acts on a Hilbert space ${\cal H}$
and $M$ is a finite-dimensional \csis\ of ${\cal H}$.
If $\Omega\in {\cal H}$ satisfies
$u\Omega=\Omega$
for some  isometry $u$ in $\con^+$,
then $\Omega\in M$.
\end{lem}
%
%
\pr
In Lemma \ref{lem:appone},
let $a_i:=u$ for all $i\geq 1$.
By assumption, we obtain $T\Omega=\Omega$.
From this and Lemma \ref{lem:appone}(ii),
$\Omega=T\Omega\in T{\cal H}\subset M$.
\qedh

%
%
\ssft{Proofs of theorems}
\label{subsection:secondtwo}
Recall ${\cal K}(\omega)$,
$\cdim\omega$ and $\kappa(\omega)$ in
$\S$ \ref{subsubsection:firstoneone}.\\

\noindent
{\it Proof of Theorem {\rm \ref{Thm:minimal}}.}
We prove the equality $\kappa(\omega)=\cdim \omega$.
This is equivalent to the following statement:
%
%
\begin{equation}
\label{eqn:twooneone}
\mbox{$\cdim \omega\leq \cdim\omega'$
for any $\omega'\in {\cal S}(\con)$
such that $\omega'\sim \omega$.}
\end{equation}
We prove (\ref{eqn:twooneone}) as follows.
Let $({\cal H},\pi,\Omega)$ denote 
the GNS representation of $\omega$.
By the assumption of $\omega(u)=1$,
we obtain $\pi(u)\Omega=\Omega$.
Assume that a state $\omega'$ on $\con$
satisfies $\omega'\sim\omega$.
Since $\omega\sim\omega'$,
we can identify $M:={\cal K}(\omega')$ with a subspace
of ${\cal H}$.
If $\dim M=\infty$, then $\cdim \omega\leq \infty=\cdim\omega'$.
Hence (\ref{eqn:twooneone}) holds.
If $\dim M<\infty$,
then $M$ and $\Omega$ 
satisfy the assumption in Lemma \ref{lem:previous}.
Hence $\Omega\in M$.
This implies
${\cal K}(\omega)\subset M={\cal K}(\omega')$.
Hence
$\cdim \omega\leq \cdim\omega'$.
\qedh

\noindent
From the proof of Theorem \ref{Thm:minimal},
the following holds.
%
%
\begin{cor}
\label{cor:gns}
For $\omega\in {\cal S}(\con)$
with GNS representation space ${\cal H}$,
assume that $\omega$ satisfies the assumption in 
Theorem {\rm \ref{Thm:minimal}} and $\cdim\omega<\infty$.
Then ${\cal K}(\omega)$ is smallest in the sense that
any nonzero finite-dimensional  \csis\ of ${\cal H}$
contains ${\cal K}(\omega)$ as a subspace.
\end{cor}
\quad 

\noindent
{\it Proof of Theorem} \ref{Thm:inclusion}.
We prove $\kappa(\omega)=\infty$.
This is equivalent to the following statement:
%
%
\begin{equation}
\label{eqn:twoonethree}
\mbox{$\cdim \omega'=\infty$ for
any $\omega'\in {\cal S}(\con)$ such that $\omega'\sim \omega$.
}
\end{equation}

Let $({\cal H},\pi,\Omega)$ denote 
the GNS representation of $\omega$.
For $l\geq 1$,
let $a[l]:=a_1\cdots a_l$ and $v_l:=\pi(a[l])^*\Omega$.
From (\ref{eqn:aoneatwo}), $X:=\{v_l:l\geq 1\}$
is an orthonormal system in ${\cal H}$
and 
%
%
\begin{equation}
\label{eqn:piacdot}
\pi(a[l]  a[l]^*)\Omega=\Omega\quad
\mbox{ for all }l\geq 1.
\end{equation}
Since $X\subset {\cal K}(\omega)$,
$\cdim\omega=\infty$.
From (\ref{eqn:piacdot}),
we obtain $T\Omega=\Omega$
where $T:=\bigwedge_{l\geq 1} \pi(a[l] a[l]^*)$.

Assume that $\omega'\in {\cal S}(\con)$
satisfies $\omega'\sim\omega$.
We identify $M:={\cal K}(\omega')$
with a subspace of ${\cal H}$.
If $\cdim \omega'<\infty$,
then $M$ satisfies 
assumptions in Lemma \ref{lem:appone}(ii).
Hence 
$\Omega=T\Omega\in T{\cal H}\subset M$.
From this,
${\cal K}(\omega)\subset M={\cal K}(\omega')$.
Hence
$\infty=\cdim \omega\leq \cdim \omega'<\infty$.
This is a contradiction.
Hence $\cdim\omega'=\infty$.
Therefore (\ref{eqn:twoonethree}) is proved.
\qedh

%
%
\sftt{Examples}
\label{section:third}
In this section, we show examples of minimal states.
For this purpose, we review properties of known states.
%
%
\ssft{Minimal states}
\label{subsection:thirdone}
%
%
\sssft{Extensions of Cuntz states}
\label{subsubsection:thirdoneone}
Recall from $\S$ \ref{subsubsection:firsttwoone} 
the definition of a Cuntz state.
For a unital C$^*$-algebra ${\cal A}$,
a unital C$^*$-subalgebra ${\cal B}$ of ${\cal A}$, 
and a state $\omega$ on ${\cal B}$,
$\omega'$ is an {\it extension} of $\omega$ to ${\cal A}$
if $\omega'$ is a state on ${\cal A}$ which satisfies
$\omega'|_{{\cal B}}=\omega$.
The following is a corollary of Theorem \ref{Thm:minimal}.
%
%
\begin{cor}
\label{cor:extensions}
Assume $2\leq n\leq m\leq \infty$.
Let $t_1,\ldots,t_m$ and $s_1\ldots,s_n$ denote
Cuntz generators of $\co{m}$ and $\con$, respectively.
Assume that there exists a unital embedding $f$ of $\co{m}$ 
into $\con$ (this requires $n\leq m$).
We identify $\co{m}$  with $f(\co{m})\subset \con$.
If $f(t_i)\in \con^+$ for all $i$,
then any extension of a Cuntz state on $\co{m}$
to $\con$ is minimal.
\end{cor}
%
%
\pr
Let $\omega$ be the Cuntz state on $\co{m}$
by $z=(z_1,\ldots,z_m)
\in ({\Bbb C}^m)_1$
and assume that
$\omega'$ is an extension of $\omega$ to $\con$.
Let $t(z):=z_1t_1+\cdots+z_mt_m$ and $u:=f(t(z))$.
Then $u\in \con^+$ and  $u^*u=I$.
By assumption,
$\omega'(u)=\omega(t(z))=1$.
From Theorem \ref{Thm:minimal},
the statement holds.
When $m=\infty$,
replace ${\Bbb C}^m$ with $\ltn$.
Then the statement holds in a similar fashion.
\qedh

%
%
\begin{prop}
\label{prop:dhjb}
There exists a representation $({\cal H},\pi)$
of $\con$ with a closed \csis\ $M$ of ${\cal H}$
such that 
$\pi(s_i){\cal H}_0$ is not a subspace of ${\cal H}_{1}$ 
for some $i\in \{1,\ldots,n\}$
where ${\cal H}=\bigoplus_{l\geq 0}{\cal H}_l$
denotes the DHJ decomposition by $M$.
\end{prop}
%
%
\pr
Let $\omega$ be the Cuntz state on $\con$
by $z=(1,0,\ldots,0)$ in $\S$ \ref{subsubsection:firstonetwo}.
Then ${\cal K}(\omega)$ equals ${\Bbb C}\Omega$
for the GNS representation $({\cal H},\pi,\Omega)$ of $\omega$,
and it is a closed \csis\ of ${\cal H}$.
For the DHJ decomposition by ${\cal K}(\omega)$,
${\cal H}_0={\Bbb C}\Omega$.
Since $\omega(s_1)=1$,
we see $\pi(s_1)\Omega=\Omega$.
This implies
$\pi(s_1){\cal H}_0={\cal H}_0\not \subset {\cal H}_1$.
\qedh

%
%
\sssft{Sub-Cuntz states}
\label{subsubsection:thirdonetwo}
Sub-Cuntz states were introduced by 
Bratteli and Jorgensen (\cite{BJ})
as extensions of Cuntz states.
We review  results in \cite{GP08}.
For $1\leq m<\infty$,
let ${\cal V}_{n,m}$ denote
the complex Hilbert space with the orthonormal basis
$\{e_J: J\in\{1,\ldots,n\}^m\}$,
that is, ${\cal V}_{n,m}
=\ell^2(\{1,\ldots,n\}^m)\cong {\Bbb C}^{n^m}$.
Let  $({\cal V}_{n,m})_1
:=\{z\in {\cal V}_{n,m}:\|z\|=1\}$.
When $n=\infty$,
let ${\cal V}_{\infty,m}:=\ell^2(\{1,2,\ldots\}^m)$.
%
%
\begin{defi}
\label{defi:subcuntz}
For $z=\sum z_Je_J\in ({\cal V}_{n,m})_1$,
$\omega$ is a {\it sub-Cuntz state} on $\con$ by $z$ 
if $\omega$ is a state on $\con$ 
which satisfies the following equations:
%
%
\begin{equation}
\label{eqn:subeqn}
\omega(s_{J})=\overline{z_{J}}
\quad \mbox{for all }J\in \{1,\ldots,n\}^m
\end{equation}
where $s_J:=s_{j_1}\cdots s_{j_m}$ when
$J=(j_1,\ldots,j_m)$,
and $\overline{z_J}$ denotes the complex conjugate of $z_J$.
In this case,
$\omega$ is called a {\it sub-Cuntz state of order} $m$.
\end{defi}

\noindent
When $n=\infty$,
replace ${\cal V}_{n,m}$ with ${\cal V}_{\infty,m}$.
A sub-Cuntz state $\omega$ of order $1$ is just a Cuntz state.

We identify 
${\cal V}_{n,m}$ with $({\cal V}_{n,1})^{\otimes m}$
by the correspondence between bases
$e_{J}\mapsto e_{j_1}\otimes \cdots \otimes e_{j_m}$
for $J=(j_1,\ldots,j_m)\in\{1,\ldots,n\}^m$.
From this identification,
we obtain
${\cal V}_{n,m}\otimes {\cal V}_{n,l}={\cal V}_{n,m+l}$
 for any $m,l\geq 1$.
Then the following hold.
%
%
\begin{Thm}
\label{Thm:maintwo}
\begin{enumerate}
\item
{\rm (\cite{GP08}, Fact 1.3)}
For any $z\in ({\cal V}_{n,m})_1$,
a sub-Cuntz state on $\con$ by $z$ exists.
\item
{\rm (\cite{GP08}, Theorem 1.4)}
For a sub-Cuntz state $\omega$ on $\con$ by $z\in ({\cal V}_{n,m})_1$,
$\omega$ is unique if and only if  $z$ is {\it nonperiodic},
that is, $z=x^{\otimes p}$ for some $x$ implies $p=1$.
In this case, $\omega$ is pure and we write it as $\tilde{\omega}_z$.
\item
{\rm (\cite{GP08}, Theorem 1.5)}
Let $p\geq 2$ and $z=x^{\otimes p}$
for a nonperiodic element $x\in ({\cal V}_{n,m'})_1$.
If $\omega$ is a sub-Cuntz state on $\con$ by $z$,
then $\omega$ is a convex hull of 
sub-Cuntz states by $e^{2\pi j\sqrt{-1}/p}x$ for $j=1,\ldots,p$.
\item
{\rm (\cite{GP08}, Theorem 1.7)}
For $z,y\in \bigcup_{m\geq 1}({\cal V}_{n,m})_1$,
assume that both $z$ and $y$ are nonperiodic.
Then the following are equivalent:
\begin{enumerate}
\item
GNS representations of $\tilde{\omega}_z$ and $\tilde{\omega}_y$
are unitarily equivalent. 
\item
$z$ and $y$ are {\it conjugate},
that is, 
$z=y$, or 
$z=x_1\otimes x_2$ and $y=x_2\otimes x_1$
for some $x_1, x_2\in \bigcup_{m\geq 1}({\cal V}_{n,m})_1$.
\end{enumerate}
\end{enumerate}
\end{Thm}

\noindent
When $n<\infty$,
any sub-Cuntz state on $\con$ is finitely correlated
(\cite{GP08}, Lemma 2.4(i)).
Furthermore, the following holds.
%
%
\begin{prop}
\label{prop:minimalsub}
Any sub-Cuntz state is minimal.
\end{prop}
%
%
\pr
For $z=\sum_{J}z_Je_J\in ({\cal V}_{n,m})_1$, let 
$\omega$ be a sub-Cuntz state on $\con$ by $z$
and let $u:=\sum_{J}z_Js_J\in \con^+$.
Then $u^*u=I$.
From (\ref{eqn:subeqn}), $\omega(u)=1$.
By Theorem \ref{Thm:minimal},
$\omega$ is minimal.
When $n=\infty$, replace ${\cal V}_{n,m}$ by ${\cal V}_{\infty,m}$. 
Then the statement is verified in a similar way.
\qedh

\noindent
Proposition \ref{prop:minimalsub}
can be also proved by using Corollary \ref{cor:extensions}
(\cite{GP15}, $\S$ 2.2).

%
%
\begin{prop}
\label{prop:differenceb}
\begin{enumerate}
\item
There exists a pure state on $\con$ which is not minimal.
\item
There exists a minimal state on $\con$ which is not pure.
\end{enumerate}
\end{prop}
%
%
\pr
(i)
Let $\omega$ be the Cuntz state by $z=(1,0,\ldots,0)$
with GNS representation $({\cal H},\pi,\Omega)$.
Let $\omega':=\omega(s_2^*(\cdot)s_2)$.
We identify ${\cal K}(\omega')$ with a subspace of ${\cal H}$.
Then ${\cal K}(\omega')$ is spanned by
the orthonormal set $\{\Omega,\pi(s_2)\Omega\}$.
Hence $\cdim\omega'=2$.
Since $\omega'\sim \omega$ and $\cdim \omega=1$,
$\omega'$ is pure but not minimal.

\noindent
(ii)
Let $\omega_{\pm}$ denote
the the Cuntz state by
$(\pm 1,0,\ldots,0)\in ({\Bbb C}^n)_1$,
respectively.
Define $\omega'':=(\omega_+ + \omega_-)/2$.
Since $\omega_+\not\sim \omega_-$
from Theorem \ref{Thm:Cuntzstateeq},
$\omega''$ is not pure.
On the other hand,
we can prove $\omega''(s_1^2)=1$.
Hence $\omega''$ is a sub-Cuntz state on $\con$
by $z=e_1^{\otimes 2}\in ({\cal V}_{n,2})_1$.
Therefore
it is minimal from Proposition \ref{prop:minimalsub}.
\qedh

%
%
\begin{prop}
\label{prop:uniqueness}
A minimal model of a state is not unique in general.
\end{prop}
%
%
\pr
Let $\omega$ and $\omega'$ be states on $\con$
such that $\omega(s_1s_2)=1=\omega'(s_2s_1)$.
Then they are pure sub-Cuntz states which exist uniquely,
and $\omega\sim \omega'$
from Theorem \ref{Thm:maintwo}(i)$\sim$(iv).
From Proposition \ref{prop:minimalsub}, they are minimal.
From $\omega'(s_2s_1)=1$,
we can prove 
$\omega'((s_2s_1)^*x)=\omega'(x)$ for any $x\in \con$.
Hence 
$\omega'(s_1s_2)=\omega'((s_2s_1)^*s_1s_2)=0$.
Therefore $\omega\ne \omega'$.
\qedh

\noindent
If $\kappa(\omega)=1$,
then a minimal model of $\omega$ is unique
from Fact \ref{fact:seven} and Theorem \ref{Thm:Cuntzstateeq}.
%
%
\begin{ex}
\label{ex:first}
{\rm
Let $\omega$ and $\omega'$ be states on $\con$
which satisfy
$\omega(s_1s_2)=1$ and $\omega'(s_1s_1+s_1s_2)=\sqrt{2}$.
Then such states are pure sub-Cuntz states
from Theorem \ref{Thm:maintwo}(ii).
From Proposition \ref{prop:minimalsub}, 
they are minimal and 
we can prove $\kappa(\omega)=\kappa(\omega')=2$, but
$\omega'\not \sim \omega\circ \alpha_g$ for any $g\in U(n)$
(see also Theorem 4.1(iv) in \cite{GP08}).
}
\end{ex}

%
%
\sssft{Geometric progression states}
\label{subsubsection:thirdonethree}
Geometric progression states were introduced 
in \cite{GP15} as extensions of Cunts states
with respect to different embeddings of Cunt algebras
from the case of sub-Cuntz states.
Assume $2\leq n<\infty$ in this section.
%
%
\begin{defi}
\label{defi:geodef}
Let $\omega\in {\cal S}(\con)$ and $m:=(n-1)k+1$ for $k\geq 2$.
\begin{enumerate}
\item
$\omega$ is 
a \gprs\ by 
$z=(z_1,\ldots,z_m)\in ({\Bbb C}^m)_1
:=\{y\in {\Bbb C}^m:\|y\|=1\}$
if $\omega$ satisfies
%
%
\begin{equation}
\label{eqn:fsubl}
\left\{
\begin{array}{rl}
\omega(s_n^rs_i)=&\overline{z_{(n-1)r+i}}\quad 
(r=0,1,\ldots,k-1, i=1,\ldots,n-1),\\
\\
\omega(s_n^k)=&\overline{z_m}.
\end{array}
\right.
\end{equation}
\item
$\omega$ is 
a \gprs\ by 
$z=(z_1,z_2,\ldots)\in\ltno:=\{y\in\ltn:\|y\|=1\}$
if $\omega$ satisfies
%
%
\begin{equation}
\label{eqn:fsubli}
\omega(s_n^{r}s_i)
=\overline{z_{(n-1)r+i}}\quad (r\geq 0,\,
i=1,\ldots,n-1).
\end{equation}
\end{enumerate}
\end{defi}
%
%
\begin{Thm}{\rm (\cite{GP15})}
\label{Thm:mainb}
\begin{enumerate}
\item
For $k\geq 2$, let $m=(n-1)k+1$.
For $z=(z_1,\ldots,z_m)\in ({\Bbb C}^m)_1$,
a \gprs\ on  $\con$ by $z$ is unique 
if and only if $|z_m|<1$.
In this case, it is pure.
We write this $\omega_z'$.
\item
For any $z\in \ltno$,
a \gprs\ on $\con$ by $z$
is unique and pure.
We write this $\omega_z'$.
\end{enumerate}
\end{Thm}
%

\noindent
In Theorem \ref{Thm:mainb}(i),
if $k=1$, then $m=n$  and $\omega$ is just a Cuntz state.
%
%
\begin{Thm}
\label{Thm:mainthree}
{\rm (\cite{GP15}, Theorem 1.8)}
Let $\omega_z'$ be as in Theorem \ref{Thm:mainb}.
\begin{enumerate}
\item
For $m=(n-1)k+1$ with $k\geq 2$,  
let ${\cal W}_m:=\{(w_1,\ldots,w_m)
\in ({\Bbb C}^m)_1:|w_m|<1\}$.
For $z,y\in {\cal W}_m$,
$\omega_z'\sim \omega_y'$ if and only if  $z=y$.
\item
For $z,y\in \ltno$,
$\omega_z'\sim \omega'_y$
if and only if $z=y$.
\end{enumerate}
\end{Thm}

%
%
\begin{Thm}
\label{Thm:reductionequiv}
{\rm (\cite{GP15}, Theorem 1.9(i))}
Let $\{e_i\}$ denote the standard basis 
of $\ltn$ and let $z\in \ltno$.
For $y=(y_1,\ldots,y_n)\in ({\Bbb C}^n)_1$,
let $\omega_y$ be as in (\ref{eqn:first}).
Then 
$\omega_z'\sim \omega_y$ if and only if 
$|y_n|<1$
and $z=\tilde{y}$
where $\tilde{y}\in \ltno$ is defined as 
%
%
\begin{equation}
\label{eqn:tilde}
\tilde{y}:=\sum_{r\geq 0}
\sum_{i=1}^{n-1}y_{n}^r y_{i}\,e_{(n-1)r+i}.
\end{equation}
\end{Thm}

%
%
\begin{Thm}
\label{Thm:mequivalence}
{\rm (\cite{GP15}, Theorem 1.10(iv))}
Assume $m=(n-1)k+1$ and $k\geq 2$.
Let $z\in {\cal W}_m$
and $y=(y_1,\ldots,y_n)\in ({\Bbb C}^n)_1$.
Let $\omega_y$ be as in (\ref{eqn:first}).
Then $\omega_z'\sim \omega_y$ 
if and only if
$|y_n|<1$ and $z=\hat{y}$
where $\hat{y}\in ({\Bbb C}^m)_1$ is defined as 
%
%
\begin{equation}
\label{eqn:tildey}
\hat{y}:=\sum_{r=0}^{k-1}
\sum_{j=1}^{n-1}y_n^{r}y_je_{(n-1)r+j}+y^k_ne_m
\end{equation}
where $\{e_j\}$ denotes the standard basis of ${\Bbb C}^m$.
\end{Thm}

%
%
\begin{Thm}
\label{Thm:finite}
{\rm (\cite{GP15}, Theorem 1.11)}
For $n<\infty$,
any geometric progression state $\omega$ 
on $\con$ of order $k<\infty$ satisfies
$\dim {\cal K}(\omega)\leq k$.
Especially,  $\omega$ is finitely correlated.
\end{Thm}

%
%
\begin{cor}
\label{cor:corgeo}
\begin{enumerate}
\item
For $z\in \ltno$,
$\kappa(\omega_z')\geq 2$
if and only if $z$ can not
be written as $\tilde{y}$ in (\ref{eqn:tilde}).
\item
For any $z\in {\cal W}_{(n-1)k+1}'$,
$\kappa(\omega_z')\leq k$.
In addition,
$2\leq \kappa(\omega_z')\leq k$
if and only if 
$z$ can not be written as $\hat{y}$ in 
(\ref{eqn:tildey}).
\end{enumerate}
\end{cor}
%
%
\pr
(i)
From Theorem \ref{Thm:reductionequiv} and 
Fact \ref{fact:seven}(ii), the statement holds.

\noindent
(ii)
Recall that 
$\kappa(\omega)\leq \cdim \omega=\dim {\cal K}(\omega)$
for any $\omega\in {\cal S}(\con)$.
From Theorem \ref{Thm:finite}, the former statement holds.
From Theorem \ref{Thm:mequivalence},
Fact \ref{fact:seven}(ii),  and the former, the latter holds.
\qedh

%
%
\begin{prop}
\label{prop:geomin}
Any \gprs\ is minimal.
\end{prop}
%
%
\pr
Assume $m=(n-1)k+1$.
For $z=(z_1,\ldots,z_m)\in ({\Bbb C}^m)_1$,
let $\omega$ be a \gprs\ on $\con$ by $z$.
Let 
$u:=\sum_{r=0}^{k-1}
\sum_{i=1}^{n-1}z_{(n-1)r+i}\,s_n^{r}s_i
+z_ms_n^k\in \con^+$.
Then $u^*u=I$.
By (\ref{eqn:fsubl}) and (\ref{eqn:fsubli}),
$\omega(u)=1$.
From Theorem \ref{Thm:minimal},
$\omega$ is minimal.

When $m=\infty$,
let 
$u:=\sum_{r\geq 0}
\sum_{i=1}^{n-1}z_{(n-1)r+i}\,s_n^{r}s_i\in \con^+$.
Then the statement holds as in the previous case. 
\qedh

\noindent
Proposition \ref{prop:geomin}
can be also proved by using Corollary \ref{cor:extensions}
(\cite{GP15}, $\S$ 1.2.3).

%
%
\ssft{Properly infinitely correlated states}
\label{subsection:thirdtwo}
In this subsection, 
we show examples of properly infinitely correlated states.
Let ${\Bbb N}:=\{1,2,\ldots\}$.
%
%
\sssft{States associated with permutative representations}
\label{subsubsection:thirdtwoone}
Let $\{e_{k,m}:(k,m)\in {\Bbb N}\times {\Bbb Z}\}$
denote the standard basis of $\ell^2({\Bbb N}\times {\Bbb Z})$.
For $2\leq n<\infty$,
define the representation $\pi$ of $\con$ 
on $\ell^2({\Bbb N}\times {\Bbb Z})$ by
%
%
\begin{equation}
\label{eqn:pisie}
\pi(s_i)e_{k,m}:=e_{n(k-1)+i,m+1}
\quad((k,m)\in {\Bbb N}\times {\Bbb Z},\,i=1,\ldots,n).
\end{equation}
By definition,
$\pi(s_1^m)^*e_{1,0}=e_{1,-m}$ for any $m\geq 1$.
Define $\omega:=\langle e_{1,0}|\pi(\cdot)e_{1,0}\rangle$.
Let $a_i:=s_1\in\con^+$ for all $i\geq 1$.
Then
%
%
\begin{equation}
\label{eqn:akalstar}
\omega(a[k]a[l]^*)=
\omega(s_1^k(s_1^*)^l)=\ip{e_{1,-k}}{e_{1,-l}}
=\delta_{k,l}\quad(l,k\geq 1).
\end{equation}
From Theorem \ref{Thm:inclusion},
$\omega$ is properly infinitely correlated.

%
%
\sssft{Induced product states}
\label{subsubsection:thirdtwotwo}
In this subsection, we review 
induced product representations \cite{ABC,BC,BC2}
and introduce induced product states.
We give a parametrization of  induced product states 
by one-sided infinite sequences of unit complex vectors. 

Let $({\Bbb C}^n)_1:=\{z\in {\Bbb C}^n:\|z\|=1\}$.
For a sequence $z\in ({\Bbb C}^n)_1^{\infty}:=
\{(z^{(i)})_{i\geq 1}:z^{(i)}\in  ({\Bbb C}^n)_1\mbox{ for all }i\}$
and $J=(j_{1},\ldots,j_{m})\in\{1,\ldots,n\}^{m}$,
define $z_J:= z^{(1)}_{j_{1}}\cdots z^{(m)}_{j_{m}}$
for $m\geq 1$
and $z_{\emptyset}:= 1$.
%
%
\begin{defi}
\label{defi:indp}
For $z=(z^{(l)})\in({\Bbb C}^n)_1^{\infty}$,
define the state $\omega_{z}$ on $\con$ as 
%
%
\begin{equation}
\label{eqn:chainstate}
\omega_{z}(s_{J}s_{K}^{*}):=
\left\{
\begin{array}{ll}
\overline{z_J}\, z_K\quad &(\mbox{when }|J|=|K|),\\
\\
0\quad &(\mbox{otherwise})
\end{array}
\right.
\end{equation}
for $J,K\in {\cal I}_n$
where $|J|$ denotes the length of a word $J$.
We call $\omega_z$ the {\it induced product state by $z$.}
\end{defi}

\noindent
Definition \ref{defi:indp} is equivalent to the original (\cite{BC},
Definition 2.9).
The GNS representation of
the state $\omega_{z}$
is called the {\it induced product representation by $z$}.
%
%
\begin{Thm}
\label{Thm:ipronetwo}
\begin{enumerate}
\item
For any $z\in({\Bbb C}^n)_1^{\infty}$, 
$\omega_z$ exists uniquely.
\item
For $z,y\in({\Bbb C}^n)_1^{\infty}$, 
$\omega_{z}\sim \omega_{y}$ if and only if 
there exists $k\geq 0$ such that 
$\sum_{l=1}^{\infty}(1-|\langle z^{(l)}|y^{(l+k)}\rangle |)<\infty$ or
$\sum_{l=1}^{\infty}(1-|\langle z^{(l+k)}|y^{(l)}\rangle |)<\infty$.
\item
For $z\in({\Bbb C}^n)_1^{\infty}$, 
$\omega_{z}$ is pure if and only if 
$\sum_{l=1}^{\infty}(1-|\langle z^{(l)}|z^{(l+k)}\rangle |)
=\infty$ for any $k\geq 1$.
In this case,
$z$ is said to be {\it aperiodic} {\rm (\cite{Laca1993})}.
\end{enumerate}
\end{Thm}
%
%
\pr
(i)
Let $\gamma$ denote the $U(1)$-gauge action on $\con$,
that is, $\gamma_z(s_i):=zs_i$ for all $i=1,\ldots,n$
and $z\in U(1)$.
Let $P$ denote the conditional expectation from $\con$
to $\con^{U(1)}:=\{x\in\con:\gamma_z(x)=x
\mbox{ for all }z\in U(1)\}\cong UHF_n$.
Let $E_{JK}:=s_Js_K^*$ for $J,K\in {\cal I}_n$.
Then $\con^{U(1)}=\overline{{\rm Lin}\langle\{E_{JK}:J,K\in 
{\cal I}_n,\, |J|=|K|\}\rangle}$.
For $z=(z^{(i)})\in ({\Bbb C}^n)_1^{\infty}$,
define the state $F_z$ on $\con^{U(1)}$ by
%
%
\begin{equation}
\label{eqn:faejk}
F_z(E_{JK}):=\overline{z_J}\,z_K
\quad (J,K\in \{1,\ldots,n\}^l,\,l\geq 1).
\end{equation}
By the natural identification ${\rm Lin}\langle\{E_{JK}:J,K\in 
\{1,\ldots,n\}^l\}\rangle\cong M_{n}({\Bbb C})^{\otimes l}$,
$F_z$ is identified with the product state
$\bigotimes_{i\geq 1}\langle z^{(i)}|(\cdot)z^{(i)}\rangle$
on  $M_{n}({\Bbb C})^{\otimes \infty}$.
Then we can verify $\omega_z=F_z\circ P$.
Hence the statement holds.

\noindent
(ii) See Theorem 3.16 of \cite{BC}.\\
(iii) See Corollary 3.17 of \cite{BC}.	
\qedh

%
%
\begin{prop}
\label{prop:indpro}
Any induced product state is properly infinitely correlated.
\end{prop}
%
%
\pr
Let $z=(z^{(i)})\in ({\Bbb C}^n)^{\infty}_1$
and $z^{(i)}=(z^{(i)}_1,\ldots,z^{(i)}_n)$.
For $i\geq 1$,
let $a_i:=\sum_{j=1}^nz^{(i)}_js_j\in \con^+$.
Then $a_i^*a_i=I$ for all $i$ and
%
%
\begin{equation}
\label{eqn:omegazal}
\omega_z(a[l]a[k]^*)
=
\sum_{|J|=l,|K|=k}
z_{J}\overline{z_{K}}
\,\omega_z(s_Js_K^*)=
\delta_{l,k}
\sum_{|J|=l=|K|}
|z_{J}|^2 |z_{K}|^2=\delta_{l,k}
\end{equation}
for $l,k\geq 1$.
From Theorem \ref{Thm:inclusion},
$\omega_z$ is properly infinitely correlated.
\qedh

%
%
\ssft{Example of an essentially finitely correlated state which is 
not finitely correlated}
\label{subsection:thirdthree}
Any finitely correlated state is essentially finitely correlated,
but the converse is not true.
%
%
\begin{prop}
\label{prop:difference}
For any $2\leq n\leq \infty$,
there exists an essentially finitely correlated state on $\con$ which 
is not finitely correlated.
\end{prop}
%
%
\pr
Let $\omega$ be the Cuntz state on $\con$
such that $\omega(s_1)=1$.
Then the following state
$\omega'$ on $\con$ 
is essentially finitely correlated, but not finitely correlated:
%
%
\begin{equation}
\label{eqn:dash}
\omega'(x):=\sum_{l\geq 1}2^{-l}\omega(A_l^* xA_l )
\quad(x\in\con)
\end{equation}
where 
$A_l:=s_2^{l-1}s_1s_2^l\in \con$ for $l\geq 1$.
In order to show this,
we prove 
%
%
\begin{equation}
\label{eqn:dashtwo}
\kappa(\omega')=1 \mbox{ and }\cdim \omega'=\infty.
\end{equation}
Since $\|A_{l}\|=1$ for all $l$ and $\omega(s_1)=1$,
we see that $\omega(A_l^* (\cdot)A_l )$ is also a state on $\con$
and it is equivalent to $\omega$ for all $l$.
Therefore $\omega'$ is also equivalent to $\omega$.
Hence we obtain $\kappa(\omega')=1$ because $\kappa(\omega)=1$.
Since $A_{l'}^*A_l=\delta_{l',l} I$,
we see 
$\omega'(x)=\omega(A^* xA)$
for $x\in\con$ where $A:=\sum_{l\geq 1}2^{-l/2}A_l\in \con$.
Let $({\cal H},\pi,\Omega)$ denote
the GNS representation of $\omega$.
For $x\in \con$, we write $\pi(x)$ as $x$ for short.
Then we can write
$\omega'=\ip{A\Omega}{(\cdot)A\Omega}$ and 
%
%
\begin{equation}
\label{eqn:dashfour}
{\cal K}(\omega')=\overline{\lspan{
s_J^* A\Omega:J\in {\cal I}_n
}}.
\end{equation}
For $l\geq 1$,
define $v_l\in {\cal K}(\omega')$ by
$v_l:=(s_2^{l-1}s_1)^*A\Omega$.
Then $v_l=2^{-l/2}s_2^l\Omega\ne 0$ for all $l$,
and 
$\ip{v_{l'}}{v_l}=2^{-(l'+l)/2}\omega((s_2^{l'})^*s_2^l)
=\delta_{l',l}/2^{l}$ because
$\omega(s_2^l)=0$ for all $l\geq 1$.
Therefore $\{v_l:l\geq 1\}$
is an infinite orthogonal system in ${\cal K}(\omega')$.
This implies $\cdim\omega'=\dim {\cal K}(\omega')=\infty$.
\qedh

%
%
\ssft{Shift representation}
\label{subsection:thirdfour}
We review the shift representation of $\con$ \cite{BJ,CFR02}.
Fix $2\leq n\leq \infty$.
Define $\Lambda:=\{1,\ldots,n\}^{\infty}$ when
$2\leq n<\infty$,
and $\Lambda:=\{1,2,\ldots\}^{\infty}$ when $n=\infty$.
Let ${\cal H}:=\ell^2(\Lambda)$
and define the representation $\Pi$ of $\con$ on ${\cal H}$ by 
%
%
\begin{equation}
\label{eqn:pisiex}
\Pi(s_i)e_x:=e_{ix}\quad(i=1,\ldots,n,\,x\in \Lambda)
\end{equation}
where $\{e_x:x\in\Lambda\}$ denotes the standard basis of ${\cal H}$
and $ix$ denotes the concatenation of two words 
$i$ and $x$ \cite{Lothaire}.
The data $({\cal H},\Pi)$ is called
the {\it shift representation} of $\con$ \cite{BJ}.
Let $\sim$ denote the {\it tail equivalence} in $\Lambda$ \cite{BJ},
that is,
for $x=(x_1,x_2,\ldots)$, $y=(y_1,y_2,\ldots)\in \Lambda$,
we write $x\sim y$ if 
there exist $p,q\geq 1$ such that 
$x_{k+p}=y_{k+q}$ for all $k\geq 1$.
For $x\in \Lambda$,
$x$ is said to be {\it \evp}
if there exist $i_0,p\geq 1$ such that 
$x_{i+p}=x_{i}$ for all $i\geq i_0$.
Otherwise, 
$x$ is said to be {\it non-\evp}.
Define $\hat{\Lambda}:=\Lambda/\!\!\! \sim$.
For $x\in \Lambda$,
we write $[x]:=\{y\in\Lambda: y\sim x\}\in\hat{\Lambda}$.
Then the following is known. 
%
%
\begin{prop}
\label{prop:shift}
\begin{enumerate}
\item
The following irreducible decomposition holds:
%
%
\begin{equation}
\label{eqn:tomegatwo}
{\cal H}
=\bigoplus_{[x]\in \hat{\Lambda}}{\cal H}_{[x]}
\end{equation}
where 
${\cal H}_{[x]}$ denotes the closed subspace of 
${\cal H}$ generated by the set $\{e_{y}:y\in [x]\}$.
\item
For $x\in \Lambda$,
let $\Pi_{[x]}$ denote the subrepresentation
of $\Pi$ associated with the subspace ${\cal H}_{[x]}$.
Then $\Pi_{[x]}$ and $\Pi_{[y]}$ are unitarily equivalent 
if and only if $x\sim y$.
Especially, {\rm (\ref{eqn:tomegatwo})} is multiplicity free.
\end{enumerate}
\end{prop}
%
%
\pr
See chapter 6 of \cite{BJ} and 
Proposition 2.5 of \cite{CFR02}.
\qedh

\noindent
In addition to Proposition \ref{prop:shift},
we show the following.
%
%
\begin{prop}
\label{prop:shiftb}
\begin{enumerate}
\item
For $x\in \Lambda$,
define $\omega_x:=\ip{e_x}{\Pi(\cdot)e_x}$
{\rm (}see also {\rm \cite{Cuntz1980}, 3.1 Proposition)}.	
Then $\omega_x$ is a pure state on $\con$ and the following hold:
\begin{enumerate}
\item
If $x\in \Lambda$ is non-\evp,  then $\omega_x$
is properly infinitely correlated, that is, $\kappa(\omega_x)=\infty$.
\item
If $x\in\Lambda$ 
is \evp\ with the period length $d$,
then $\kappa(\omega_{x})=d$.
\end{enumerate}
\item
Let $\omega_x$ be as in {\rm (i)}.
If $x\in\Lambda$ is \evp, then the following are equivalent:
\begin{enumerate}
\item
$\omega_x$ is minimal.
\item
$x=(x_1,x_2,\ldots)$ is purely periodic,
that is, 
there exists $d\geq 1$ such that
$x_{i+d}=x_i$ for all $i\geq 1$.
\end{enumerate}
\end{enumerate}
\end{prop}
%
%
\pr
(i)
The purity of $\omega_x$ holds from the irreducibility 
of $\Pi_{[x]}$.\\
(a) 
Assume that $x=(x_1,x_2,\ldots)\in \Lambda$ is non-\evp.
Define $a_i:=s_{x_i}\in \con^+$.
Then $a_i^*a_i=I$ for all $i$.
Since $x$ is non-\evp, 
we see that $\omega_x(a[k]a[l]^*)=\delta_{kl}$
for all $k,l\geq 1$.
From Theorem \ref{Thm:inclusion},
$\omega_x$ is properly infinitely correlated.

\noindent
(b)
Assume that $x\in\Lambda$ has a minimal repeating block 
$x^{'}\in \{1,\ldots,n\}^d$.
Then there exists $x''\in  \{1,\ldots,n\}^c$
such that $x=x''x'x'x'\cdots$.
Remark that $x'$ is not periodic by definition.
Define $\hat{x}:=x'x'x'\cdots\in\Lambda$.
Then $\omega_x\sim \omega_{\hat{x}}$ because $x\sim \hat{x}$
and Proposition \ref{prop:shift}(ii).
Let $({\cal H}',\pi',\Omega)$ denote
the GNS representation of $\omega_{\hat{x}}$.
When $x'=(j_1,\ldots,j_d)$,
let $v_k:=\pi'(s_{j_k}\cdots s_{j_d})\Omega$ for $k=1,\ldots,d$.
Then we can verify that 
$\{v_k:k=1,\ldots,d\}$ is an orthonormal basis of 
${\cal K}(\omega_{\hat{x}})$
because $x'$ is not periodic.
Therefore $\cdim \omega_{\hat{x}}=d$.
Let $u:=s_{j_1}\cdots s_{j_d}\in \con^+$.
Then $u^*u=I$ and $\omega_{\hat{x}}(u)=1$.
Therefore $\omega_{\hat{x}}$ is minimal
from Theorem \ref{Thm:minimal}.
From this,
$\kappa(\omega_x)
=\kappa(\omega_{\hat{x}})
=\cdim \omega_{\hat{x}}=d$.

\noindent
(ii)
Let  $d$ be the period length of $x$.
Let $x=(x_1,x_2,\ldots)$
and define
$x^{(i)}:=(x_i,x_{i+1},\ldots)\in\Lambda$ for $i\geq 1$.
From (i)(b), $\kappa(\omega_x)=d$
and $x^{(i+d)}=x^{(i)}$ for $i\geq i_0$ for some $i_0\geq 1$.
From Proposition \ref{prop:shift}(i) and (ii),
we can identify ${\cal K}(\omega_x)$
with a subspace of ${\cal H}_{[x]}$
generated by $X_x:=\{\Pi(s_J)^*e_x:J\in {\cal I}_n\}
\setminus \{0\}$.
From (\ref{eqn:pisiex}), 
we see $X_x=\{e_{x^{(i)}}:i\geq 1\}$.
From this and 
$\ip{e_{x^{(i)}}}{ e_{x^{(j)}}}=\delta_{x^{(i)},x^{(j)}}$,
we obtain
%
%
\begin{equation}
\label{eqn:xxx}
\cdim \omega_x=\#X_x.
\end{equation}

\noindent
(a)$\Rightarrow$(b).
Assume that $\omega_x$ is minimal.
From this and (\ref{eqn:xxx}),
$\#X_x=\cdim \omega_x=\kappa(\omega_x)=d$.
Therefore $X_x=\{e_{x^{(1)}},\ldots,e_{x^{(d)}}\}$.
Hence $x$ is purely periodic.

\noindent
(b)$\Rightarrow$(a).
Assume that $x$ is purely periodic.
Then $X_x=\{e_{x^{(1)}},\ldots,e_{x^{(d)}}\}$.
From (\ref{eqn:xxx}),
$\cdim \omega_x=\#X_x=d=\kappa(\omega_x)$.
Hence $\omega_x$ is minimal.
\qedh

%
%
\begin{rem}
\label{rem:dynamical}
{\rm
We give an interpretation of the invariant $\kappa$ 
as the theory of symbolic dynamical systems
from Proposition \ref{prop:shiftb}.
For an \evp\ element $x\in \Lambda$,
let $d(x)$ denote the {\it period length} of $x$, that is,
the length of a minimal repeating block of $x$.
For a non-\evp\ element $x\in \Lambda$,
we define $d(x):=\infty$.
Then the map 
%
%
\begin{equation}
\label{eqn:invd}
d:\Lambda\to \{1,2,\ldots,\infty\}
\end{equation}
is surjective, and if $x\sim y$, then $d(x)=d(y)$,
that is,  $d$ is an invariant of elements in the orbit space 
$\hat{\Lambda}$.
By using $\kappa$,
we can write
%
%
\begin{equation}
\label{eqn:kappadyna}
\kappa(\omega_x)=d(x)\quad(x\in \Lambda)
\end{equation}
where $\omega_x$ denotes the state 
in Proposition \ref{prop:shiftb}(i).
Therefore the invariant $\kappa(\omega)$ 
of a state $\omega$
can be regarded as a generalization of the period length of 
an orbit of the full one-sided shift on $\Lambda$ \cite{Kitchens}.
This perspective is natural in a sense that
a Cuntz algebra is a special Cuntz-Krieger algebra \cite{CK}
and
Cuntz-Krieger algebras were introduced as 
a class of C$^{*}$-algebra associated with topological Markov chains.
From Theorem \ref{prop:shiftb}(ii),
the minimality of a state is also interpreted as the pure periodicity 
of an element in $\Lambda$.
}
\end{rem}

%
%
\ssft{Cardinality of minimal pure states}
\label{subsection:thirdfive}
%
%
\begin{Thm}
\label{Thm:conti}
For any $2\leq n\leq \infty$ and $1\leq d\leq \infty$,
there exist continuously many mutually inequivalent 
pure states $\omega$ on $\con$ 
which are minimal and $\kappa(\omega)=d$.
\end{Thm}
We prove Theorem \ref{Thm:conti} as follows.
%
%
\sssft{$d<\infty$}
\label{subsubsection:thirdfiveone}
Fix $1\leq d<\infty$.

Assume $n<\infty$.
Let ${\cal V}_{n,m}$ be as in $\S$ \ref{subsubsection:thirdonetwo}.
We identify 
${\cal V}_{n,m}$ with $({\cal V}_{n,1})^{\otimes m}
=\lspan{e_{i_1}\otimes \cdots \otimes e_{i_m}:i_1,\ldots,i_{m}=1,\ldots,n}$.
For $c\in U(1):=\{c\in {\Bbb C}:|c|=1\}$, 
let $\rho_c$ denote the sub-Cuntz state on $\con$ by
$z=\overline{c}\, e_2^{\otimes (d-1)}\otimes e_1\in
({\cal V}_{n,d})_1$.
From Theorem \ref{Thm:maintwo}(ii),
$\rho_c$ is uniquely defined as a state which satisfies
%
%
\begin{equation}
\label{eqn:stwo}
\rho_c(s_2^{d-1}s_1)=c,
\end{equation}
and it is pure.
Let $({\cal H},\pi,\Omega)$ denote
the GNS representation of $\rho_c$.
For $x\in \con$, we write $\pi(x)$ as $x$ for short.
From (\ref{eqn:stwo}), we obtain 
$s_2^{d-1}s_1\Omega=c\Omega$.
Let $v_{i}:=s_2^{i-1}s_1\Omega\in {\cal H}$
for $i=1,\ldots,d$.
Then we can verify that $\{v_1,\ldots,v_d\}$ 
is an orthonormal basis of ${\cal K}(\rho_c)$.
Hence we obtain $\cdim\rho_c=d$.
From Proposition \ref{prop:minimalsub},
$\kappa(\rho_c)=\cdim\rho_c=d$.
From Theorem \ref{Thm:maintwo}(iv),
$\rho_c\sim \rho_{c'}$ if and only if $c=c'$.

When $n=\infty$, replace ${\cal V}_{n,m}$ 
with ${\cal V}_{\infty,m}$.
Then the statement is verified in a similar way.

Hence Theorem \ref{Thm:conti} holds
when $d<\infty$.

%
%
\sssft{$d=\infty$}
\label{subsubsection:thirdfivetwo}
Let $\Lambda,\sim, \hat{\Lambda},[x]$
be as in $\S$ \ref{subsection:thirdfour}.
We write 
$\aleph_0$ and $\aleph_1$ for 
the cardinalities of ${\Bbb N}$ and ${\Bbb R}$,
respectively.
%
%
\begin{lem}
\label{lem:continuous}
\begin{enumerate}
\item
For any $x\in \Lambda$,
$\#[x]=\aleph_0$.
\item
$\#\hat{\Lambda}=\aleph_1$.
\item
Let  
$\hat{\Lambda}_{{\rm nep}}
:=\{[x]\in \hat{\Lambda}:
x \mbox{ is non-\evp}\}$.
Then
$\#\hat{\Lambda}_{{\rm nep}}=\aleph_1$.
\end{enumerate}
\end{lem}
%
%
\pr
Define
$A^+:=\bigcup_{l\geq 1}\{1,\ldots,n\}^l$ \cite{Lothaire}.
When $n=\infty$,
replace $\{1,\ldots,n\}^l$ with  $\{1,2,\ldots\}^{l}$ 
for each $l\geq 1$.

\noindent
(i)
For $x\in \Lambda$,
if $y\in [x]$,
then $y=y_1x_2$ for some $y_1,x_1\in A^+$
and $x_2\in \Lambda$ such that $x=x_1x_2$.
From this, $y$ is determined  only by $y_1$.
Hence $[x]\cong A^+$ as a set.
Therefore $\#[x]=\#A^+ =\aleph_0$.

\noindent
(ii)
Since $\#\Lambda=\aleph_1$,
the statement holds from (i).

\noindent
(iii)
Let $\hat{\Lambda}_{{\rm ep}}:=\{[x]\in \hat{\Lambda}:
x\mbox{ is \evp}\}$.
Then $\hat{\Lambda}=\hat{\Lambda}_{{\rm ep}}
\sqcup \hat{\Lambda}_{{\rm nep}}$.
Any $[x]\in \hat{\Lambda}_{{\rm ep}}$ has
a minimal repeating block $x'\in A^+$.
Hence
$\hat{\Lambda}_{{\rm ep}}\cong A^+$ as a set.
Therefore 
$\#\hat{\Lambda}_{{\rm ep}}=\#A^+=\aleph_0$.
Hence $\#\hat{\Lambda}_{{\rm nep}}=
\#(\hat{\Lambda}\setminus \hat{\Lambda}_{{\rm ep}})
=\aleph_1$
from (ii).
\qedh

\noindent
From Proposition \ref{prop:shift}(ii) and 
Proposition \ref{prop:shiftb}(i),
$\{\omega_x:[x]\in\hat{\Lambda}\}$
is a set of mutually inequivalent pure states on $\con$.
From this and Proposition \ref{prop:shiftb}(i)(a),
$\Xi:=\{\omega_x:[x]\in\hat{\Lambda}_{{\rm nep}}\}$
is a set of mutually inequivalent properly infinitely correlated pure states 
on $\con$.
Since a properly infinitely correlated state $\omega$ satisfies
$\infty=\kappa(\omega)\leq \cdim\omega\leq \infty$,
it is minimal.
From this,
$\Xi$ is a set of mutually inequivalent minimal pure states on $\con$
with $\kappa=\infty$.
From this  and Lemma \ref{lem:continuous}(iii),
$\#\Xi=\#\hat{\Lambda}_{{\rm nep}}=\aleph_1$.
Hence the case of $d=\infty$ in Theorem \ref{Thm:conti}
is proved.
%
%
\sftt{Applications}
\label{section:fourth}
%
%
\ssft{Invariant of irreducible representations of $\con$}
\label{subsection:fourthone}
Let $\Irr\con$ denote
the class of all irreducible representations of $\con$.
For $\pi,\pi'\in \Irr\con$,
we write $\pi\sim\pi'$ when
they are unitarily equivalent.
As an application of $\kappa$ in (\ref{eqn:cdimkappa}),
we introduce an invariant of arbitrary irreducible representations
of $\con$ with respect to $\sim$.

For $({\cal H},\pi)\in \Irr\con$
and $x\in {\cal H}_1:=\{y\in {\cal H}:\|y\|=1\}$,
define
%
%
\begin{equation}
\label{eqn:pikappa}
\kappa(\pi):=\kappa(\omega_x\circ \pi)
\end{equation}
where $\omega_x:=\ip{x}{(\cdot)x}$.
Then $\kappa(\pi)$ is independent in the choice of
$x$ because $\omega_x\circ \pi\sim \omega_y\circ \pi$
for any $y\in {\cal H}_1$.
From Theorem \ref{Thm:mainone}(ii) and 
Proposition \ref{prop:unitary}(ii),
the following hold.
%
%
\begin{prop}
\label{prop:irreducible}
\begin{enumerate}
\item
For $\pi,\pi'\in \Irr\con$,
if $\pi\sim\pi'$,
then $\kappa(\pi)=\kappa(\pi')$.
\item
For any $\pi\in \Irr\con$ and $g\in U(n)$,
$\kappa(\pi\circ \alpha_g)=\kappa(\pi)$.
\end{enumerate}
\end{prop}

\noindent
By combining 
Proposition \ref{prop:irreducible}(i) and (ii),
if $\pi,\pi'\in \Irr\con$ satisfy
$\pi\sim\pi'\circ \alpha_g$ for some $g\in U(n)$,
then $\kappa(\pi)=\kappa(\pi')$.

For example,
$\Pi_{[x]}$ in Proposition \ref{prop:shift}(ii)
and
$d(x)$ in Remark \ref{rem:dynamical}
satisfy
$\kappa(\Pi_{[x]})=d(x)$ for all $x\in\Lambda$.

Let $\widehat{\con}$ denote the 
{\it spectrum} of $\con$ \cite{Ped},
that is, the set of all unitary equivalence
classes of irreducible representations of $\con$.
Then we obtain the following decomposition by using $\kappa$:
%
%
\begin{equation}
\label{eqn:widecon}
\widehat{\con}=\coprod_{d=1}^{\infty}\widehat{\con}(d),
\quad 
\widehat{\con}(d):=\{[\pi]\in \widehat{\con}:
\kappa(\pi)=d\}\quad(1\leq d\leq\infty)
\end{equation}
where $[\pi]:=\{\pi'\in \Irr\con: \pi'\sim \pi\}$.
From Theorem \ref{Thm:conti},
$\#\widehat{\con}(d)=\aleph_1$
for all $2\leq n\leq\infty$ and $1\leq d\leq\infty$.
Especially, $\widehat{\con}(1)\cong$
``the set of all Cuntz states on $\con$"
by Fact \ref{fact:seven}.

%
%
\ssft{New invariant of ergodic endomorphisms 
of ${\cal B}({\cal H})$}
\label{subsection:fourthtwo}
For ${\cal H}:=\ltn$, let $\endbh$
denote the set of all unital endomorphisms of ${\cal B}({\cal H})$.
For $\varphi_1,\varphi_2\in \endbh $,
$\varphi_1$ and $\varphi_2$ are said to be {\it conjugate}
if there exists an automorphism $\gamma$ of ${\cal B}({\cal H})$
such that $\varphi_2=\gamma\circ \varphi_1\circ \gamma^{-1}$.
In this case, we write $\varphi_1\sim\varphi_2$.
The classification problem of elements in 
$\endbh $ by $\sim$
has been considered in
\cite{Arveson1989, BJ1997,BJP,FL,Laca1993B,PR}.
As an invariant of elements in $\endbh$,
the Powers index is well known \cite{Pow1988}.
We introduce a new invariant for a special subset of $\endbh$.
%
%
\begin{Thm}
\label{Thm:powers}
{\rm (\cite{Arveson1989,Laca1993B}}, see also
 {\rm $\S$ 3} of {\rm \cite{BJP})}
For $2\leq n\leq \infty$,
let $s_1,\ldots,s_n$ denote Cuntz generators of $\con$ and 
let $\co{1}:=C({\Bbb T})$.
For $\co{1}$, we define $s_1$ as a unitary which generates $\co{1}$.
Let $\Rep(\con,{\cal H})$
denote the set of all unital representations of $\con$ on ${\cal H}$.
For any $\varphi\in \endbh $,
there exist $1\leq n\leq \infty$ 
and $\pi\in \Rep(\con,{\cal H})$
such that $\varphi=\sum_{i=1}^n\pi(s_i)(\cdot)\pi(s_i)^*$.
The number $n$ is called the {\it Powers index} of $\varphi$.
We write $\Ind\varphi$ as $n$.
\end{Thm}

Assume $2\leq n\leq \infty$.
Let 
%
%
\begin{equation}
\label{eqn:endn}
\endnbh:=
\{\varphi\in \endbh : \Ind\varphi=n\}.
\end{equation}
For $\pi \in \Rep(\con,{\cal H})$,
define $\varphi_{\pi}\in \endbh$ by
%
%
\begin{equation}
\label{eqn:varphisum}
\varphi_{\pi}:=\sum_{i=1}^n \pi(s_i)(\cdot)\pi(s_i)^*.
\end{equation}
From this and Theorem \ref{Thm:powers}, the map
%
%
\begin{equation}
\label{eqn:repcon}
\Rep(\con,{\cal H})\ni \pi\mapsto \varphi_{\pi}
\in \endnbh
\end{equation}
is surjective. 
In other words, we can write
$\endnbh=\{\varphi_{\pi}:\pi \in \Rep(\con,{\cal H})\}$.
About this map,  the following holds.
%
%
\begin{Thm}
\label{Thm:equivend}
{\rm (\cite{Arveson1989,Laca1993B}}, see also 
{\rm $\S$ 3} of {\rm \cite{BJP})}
Assume $2\leq n,m\leq \infty$.
For $\pi_1\in \Rep(\con,{\cal H})$
and $\pi_2\in \Rep(\co{m},{\cal H})$,
the following are equivalent:
\begin{enumerate}
\item
$\varphi_{\pi_1}\sim \varphi_{\pi_2}$.
\item
$n=m$ 
and $\pi_1\sim\pi_2\circ \alpha_g$ for some $g\in U(n)$.
\end{enumerate}
Especially, 
$\varphi_{\pi_1}= \varphi_{\pi_2}$
if and only if 
$n=m$ and $\pi_1=\pi_2\circ \alpha_g$ for some $g\in U(n)$.
\end{Thm}

\noindent
Remark that
$\endnbh$ is in one-to-one correspondence with
the $U(n)$-orbit space
$\Irr(\con,{\cal H})/U(n)
:=\{\langle\pi\rangle: \pi\in \Irr(\con,{\cal H})\}$
where $\langle \pi\rangle:=\{\pi\circ \alpha_g: g\in U(n)\}$.
%
%
\begin{Thm}
\label{Thm:ergodic}
{\rm (\cite{BJP})}
Assume $2\leq n\leq \infty$.
For $\pi\in \Rep(\con,{\cal H})$,
the following are equivalent:
\begin{enumerate}
\item
$\varphi_{\pi}$ is {\it ergodic}, that is,
$\{X\in {\cal B}({\cal H}): \varphi_{\pi}(X)=X\}={\Bbb C}I$.
\item
$\pi$ is irreducible.
\end{enumerate}
\end{Thm}

By combining Theorem \ref{Thm:ergodic},
Theorem \ref{Thm:equivend}, and (\ref{eqn:pikappa}),
we can define a number $\kappa(\varphi_{\pi})$ 
for an ergodic endomorphism $\varphi_{\pi}$ by
%
%
\begin{equation}
\label{eqn:kappaphi}
\kappa(\varphi_{\pi}):=\kappa(\pi)
\end{equation}
where $\kappa(\pi)$ is as in (\ref{eqn:pikappa}).
From Proposition \ref{prop:irreducible}(ii),
$\kappa(\varphi_{\pi})$ is well defined.
From Proposition \ref{prop:irreducible}(i),
$\varphi_{\pi}\sim \varphi_{\pi'}$ implies 
$\kappa(\varphi_{\pi})=\kappa(\varphi_{\pi'})$,
that is, we obtain an invariant 
of ergodic endomorphisms: 
%
%
\begin{equation}
\label{eqn:rmerg}
\kappa:\Erg_n{\cal B}({\cal H})\to \{1,2,\ldots,\infty\}
\end{equation}
where $\Erg_n{\cal B}({\cal H}):=\{
\varphi\in\endnbh:
\varphi\mbox{ is ergodic}
\}$.
From examples in $\S$ \ref{section:third},
we can construct an ergodic endomorphism $\varphi$ with 
$\Ind\varphi=n$ and $\kappa(\varphi)=d$
for any 
$2\leq n\leq\infty$ and 
$1\leq d\leq \infty$.

\ww {\bf Acknowledgment:}
The author would like to express his sincere thanks 
to Kengo Matsumoto
for a relation between \cite{Cuntz1980} and the shift representation,
and would like to thank the referees for the careful reading
of the draft.

%
%

%
%

%
\label{Lastpage}
\end{document}